\documentclass[a4paper,12pt,leqno]{amsart}
\usepackage{amsmath,amssymb,mathrsfs}
\usepackage[matrix,arrow]{xy} 
\usepackage{xcolor,hyperref}
\hypersetup{colorlinks=true,citecolor=black,linkcolor=black}
\newcommand{\h}{\hbox}
\newcommand{\q}{\quad}
\newcommand{\nin}{\noindent}

\newcommand{\ms}{\par\medskip}
\newcommand{\sk}{\par\smallskip}

\newcommand{\msn}{\par\medskip\noindent}

\newcommand{\ges}{\geqslant}
\newcommand{\les}{\leqslant}
\newcommand{\one}{\hskip1pt}

\newcommand{\msum}{\hbox{$\sum$}}

\newcommand{\mprod}{\hbox{$\prod$}}
\newcommand{\D}{{\mathcal D}}
\newcommand{\E}{{\mathcal E}}
\newcommand{\G}{{\mathcal G}}
\newcommand{\K}{{\mathcal K}}
\newcommand{\Rc}{{\mathcal R}}

\newcommand{\OO}{{\mathcal O}}

\newcommand{\Q}{{\mathbb Q}}
\newcommand{\C}{{\mathbb C}}
\newcommand{\N}{{\mathbb N}}
\newcommand{\R}{{\mathbb R}}

\newcommand{\Z}{{\mathbb Z}}

\newcommand{\Ft}{\widetilde{F}}
\newcommand{\Wt}{\widetilde{W}}
\newcommand{\Gr}{{\rm Gr}}

\newcommand{\aTj}{\alpha^{\rm\hskip-1pt Tj}}
\newcommand{\Sb}{{\bf S}}

\newcommand{\Cf}{C_{\hskip-1pt f}}
\newcommand{\If}{I_{\hskip-1pt f}}
\newcommand{\Tf}{T_{\!f}}
\newcommand{\Gf}{\G_{\hskip-1pt f}}
\newcommand{\Hsf}{H''_{\hskip-1pt f}}
\newcommand{\dti}{\dd_t^{-1}}
\newcommand{\dt}{\dd_t}
\newcommand{\al}{\alpha}
\newcommand{\be}{\beta}
\newcommand{\ga}{\gamma}
\newcommand{\Ga}{\Gamma}

\newcommand{\De}{\Delta}
\newcommand{\la}{\lambda}

\newcommand{\si}{\sigma}
\newcommand{\Si}{\Sigma}
\newcommand{\ep}{\varepsilon}
\newcommand{\om}{\omega}
\newcommand{\Om}{\Omega}
\newcommand{\dd}{\partial}
\newcommand{\ddd}{{\rm d}}
\newcommand{\Sp}{{\rm Sp}}
\newcommand{\SpT}{{\rm Sp}^{\mhs\rm Tj}}
\newcommand{\Spf}{{\rm Sp}_{\mhs f}}
\newcommand{\SpTf}{{\rm Sp}\mhs^{\rm Tj}_f}
\newcommand{\mm}{{\mathfrak m}}
\newcommand{\tos}{\,{\to}\,}
\newcommand{\eq}{\,{=}\,}
\newcommand{\defs}{\,{:=}\,}
\newcommand{\nes}{\,{\ne}\,}
\newcommand{\ins}{\,{\in}\,}
\newcommand{\sst}{\,{\subset}\,}
\newcommand{\stm}{\,{\setminus}\,}
\newcommand{\gess}{\,{\ges}\,}
\newcommand{\less}{\,{\les}\,}
\newcommand{\sgt}{\,{>}\,}
\newcommand{\slt}{\,{<}\,}
\newcommand{\col}{\,{:}\,}
\newcommand{\pl}{\one {+}\one}
\newcommand{\mi}{\one {-}\one}
\newcommand{\bl}{\bigl}
\newcommand{\br}{\bigr}

\newcommand{\ssb}{\raise.15ex\h{${\scriptscriptstyle\bullet}$}}
\newcommand{\ssc}{\,\raise.15ex\h{${\scriptstyle\circ}$}\,}
\newcommand{\onto}{\twoheadrightarrow}
\newcommand{\into}{\hookrightarrow}
\newcommand{\simto}{\,\,\rlap{\hskip1.5mm\raise1.4mm\hbox{$\sim$}}\hbox{$\longrightarrow$}\,\,}
\newcommand{\mhs}{\hskip-1pt}

\makeatletter
\renewcommand\section{\@startsection{section}{1}{0pt}{-3ex plus -1ex minus -.2ex}{2.3ex plus.2ex}{\centering\normalfont\bfseries}}
\makeatother
\theoremstyle{plain}
\newtheorem{thm}{Theorem}[section]
\newtheorem{cor}[thm]{Corollary}
\newtheorem{prop}[thm]{Proposition}

\newtheorem{ithm}{Theorem}
\newtheorem{icor}{Corollary}
\newtheorem{iprop}{Proposition}

\theoremstyle{definition}
\newtheorem{rem}[thm]{Remark}

\newtheorem{irem}{Remark}

\begin{document}
\title[Tjurina spectrum and graded symmetry]
{Tjurina spectrum and graded symmetry\\of missing spectral numbers}
\author[S.-J. Jung]{Seung-Jo Jung}
\address{S.-J. Jung : Department of Mathematics Education, and Institute of Pure and Applied Mathematics, Jeonbuk National University, Jeonju, 54896, Korea}
\email{seungjo@jbnu.ac.kr}
\author[I.-K. Kim]{In-Kyun Kim}
\address{I.-K. Kim : June E Huh Center for Mathematical Challenges, Korea Institute for Advanced Study, 85 Hoegiro Dongdaemun-gu, Seoul 02455, Korea}
\email{soulcraw@kias.re.kr}
\author[M. Saito]{Morihiko Saito}
\address{M. Saito : RIMS Kyoto University, Kyoto 606-8502 Japan}
\email{msaito@kurims.kyoto-u.ac.jp}
\author[Y. Yoon]{Youngho Yoon}
\address{Y. Yoon : Department of Mathematics, Chungbuk National University, Cheongju-si, Chungcheongbuk-do, 28644, Korea}
\email{mathyyoon@gmail.com}
\thanks{This work was partially supported by National Research Foundation of Korea (the first author: NRF-2021R1C1C1004097, the second author: NRF-2023R1A2C1003390 and NRF-2022M3C1C8094326, and the fourth author: RS-2023-00245670).}
\begin{abstract} For a hypersurface isolated singularity defined by a convergent power series $f$, the Steenbrink spectrum can be defined as the Poincar\'e polynomial of the graded quotients of the $V$-filtration on the Jacobian ring of $f$. The Tjurina subspectrum is defined by replacing the Jacobian ring with its quotient by the image of the multiplication by $f$. We prove that their difference (consisting of missing spectral numbers) has a canonical graded symmetry. This follows from the self-duality of the Jacobian ring, which is compatible with the action of $f$ as well as the $V$-filtration. It implies for instance that the number of missing spectral numbers which are smaller than $(n{+}1)/2$ (with $n$ the number of variables) is bounded by $[(\mu{-}\tau)/2]$. We can moreover improve the estimate of Brian\c{c}on-Skoda exponent in the semisimple monodromy case.
\end{abstract}
\maketitle

\section*{Introduction} \label{intr}
\nin
Tjurina number is an interesting complex analytic invariant of an isolated hypersurface singularity, which measures in some sense how far the singularity is away from being weighted homogeneous \cite{SaK}. However, it does not seem necessarily trivial to construct a good refinement related to Hodge theory.
\sk
Let $f\ins\C\{x\}$ be a convergent power series having an isolated singularity at the origin with $f(0)\eq0$. Here $x\eq(x_1,\dots,x_n)$ is the coordinate system of $(\C^n,0)$.
The {\it Steenbrink spectrum\one} $\Spf(t)\eq\msum_{i=1}^{\mu}\,t^{\one\al_i}$ (with $\mu$ the Milnor number) can be defined as the Poincar\'e polynomial of the $V$-filtration on the Jacobian ring $\C\{x\}/(\dd f)$ so that
\begin{equation} \label{1}
|\If^{\al}|=\dim_{\C}\Gr^{\al}_V\bl(\C\{x\}/(\dd f)\br)\q\h{with}\q\If^{\al}\eq\{i\in\If\mid\al_i\eq\al\}.
\end{equation}
Here $\If\defs\{1,\dots,\mu\}$, $(\dd f)\sst\C\{x\}$ is the Jacobian ideal, and $V$ is the quotient filtration of the $V$-filtration on the Brieskorn lattice (or the {\it microlocal\one} $V$-filtration on $\C\{x\}$), see \cite{SS}, \cite{Va1} (and also \cite{JKSY}). Here the spectral numbers $\al_i$ are assumed to be {\it weakly increasing.} We may assume that $f$ is a polynomial by the finite determinacy, see for instance \cite{GLS}.
\sk
We define the {\it Tjurina subspectrum\one} (or {\it spectrum\one}) $\SpTf(t)=\msum_{j=1}^{\tau}\,t^{\one\aTj_j}=\msum_{i\in\Tf}\,t^{\one\al_i}$ with $\Tf\sst\If$ by
\begin{equation*}
|\Tf^{\al}|=\dim_{\C}\,\Gr^{\al}_V\bl(\C\{x\}/(\dd f,f)\br)\q\h{with}\q\Tf^{\al}\eq\If^{\al}\cap\Tf\eq\If^{\al}\cap [1,i_{\al}],
\end{equation*}
where $i_{\al}\ins\Tf^{\al}$, and the $\aTj_j$ are weakly increasing, see \cite{JKY}, \cite{JKSY}. Note that $|\Tf|\eq\SpTf(1)$ is equal to the Tjurina number $\tau$ of $f$. Set $\Cf\defs\If\stm\Tf$, and
\begin{equation*}
\Spf^C(t)\defs\Spf(t)\mi\SpTf\eq\msum_{i\in\Cf}\,t^{\one\al_i}\eq\msum_{k=1}^{\mu-\tau}\,t^{\one\al^C_k},
\end{equation*}
where the $\al^C_k$ are assumed to be weakly increasing. This is the Poincar\'e polynomial of the image of the endomorphism
\begin{equation*}
[f]\ins{\rm End}_{\C\{x\}}\bl(\C\{x\}/(\dd f)\br),
\end{equation*}
where $[f]$ denotes the class of $f$. It is called the {\it complemental part\one} to the Tjurina spectrum. The numbers $\al^C_k$ ($k\ins[1,\mu{-}\tau]$) and $\al_i$ ($i\ins\Cf$) are called the {\it missing spectral numbers\one} of $f$. Note that some of the missing spectral numbers may also belong to the Tjurina spectrum, for instance, if $f\eq x^6\pl y^6\pl x^4y^3$, where $\If^{3/2}\eq\{23,24\}$ and $\Tf^{3/2}\eq\{23\}$.
\sk
In this paper we {\it assume}
\begin{equation*}
\mu\nes\tau.
\end{equation*}
This is equivalent to assuming that $f$ is not a weighted homogeneous polynomial \cite{SaK}. The Tjurina spectrum is a quite subtle invariant of a hypersurface isolated singularity with $\mu\nes\tau$, see Proposition\,\,\ref{P4.3} below. Recall that the hypersurface germ can be determined by using the Tjurina algebra $\C\{x\}/(\dd f,f)$, see \cite{MY}.
(The relation between the Tjurina spectrum and the Bernstein-Sato polynomial does not seem very clear, see for instance \cite{sem} and also Remark\,\,\ref{R4.3} below.)

\begin{irem}[{\it Instability under $\tau$-constant deformations\one\rm}] \label{R1}
The Tjurina spectrum is {\it unstable under a $\one\tau$-constant deformation,} for instance, if $f_u\eq x^6\pl y^5\pl u\one x^4y^2\pl x^3y^3$ ($u\ins\C$), where the missing spectral numbers for $f_0$ and $f_u$ ($u\nes0$) are respectively $\tfrac{44}{30},\tfrac{49}{30}$ and $\tfrac{43}{30},\tfrac{49}{30}$ using Remark\,\,\ref{R3.5} below. (Indeed, $f_u\eq u\one c\one x^4y^2\pl c'x^3y^3$ in $\C\{x\}/(\dd f_u$) for some $c,c'\ins\C^*$.) So it does not seem easy to determine the minimal stratification of the base space of a miniversal $\mu$-constant deformation of $f$ such that the Tjurina spectrum is constant on each stratum in general.
\end{irem}

Using a spectral sequence together with the self-duality of the Jacobian ring, which is compatible with the multiplication by $[f]$ and also with the $V$-filtration, we can prove the following.

\begin{ithm} \label{T1}
Set $A\defs\C\{x\}/(\dd f)$, $I^{\al}\!A\defs[f]V^{\al}\!A\,\,(\al\ins\Q)$. For $\al,\be,\ga\ins\Q$ with $\al{+}\be\eq\ga{+}n$, $\ga\gess1$, there is the canonical perfect pairing
\begin{equation*}
\Gr_I^{\al-\ga}\Gr_V^{\al}A\otimes_{\C}\Gr_I^{\be-\ga}\Gr_V^{\be}A\to\C.
\end{equation*}
\end{ithm}

Here the bigraded pieces $\Gr_I^{\al-\ga}\Gr_V^{\al}A$ vanish for $\ga\slt1$ by Proposition\,\,\ref{P1.1} below. (Note that there is no relation between $I^{\al}$ and $\If^{\al}$.) Theorem\,\,\ref{T1} is closely related to the {\it self-duality of spectral sequences\one} (see \cite[(1.3.1)]{De}, \cite[5.2.19--20]{mhp}) and may follow from it although a direct argument is adopted instead of it, see Remark\,\,\ref{R5.1} below. (This does not seem to be a standard technique of this field.)
In the case $\Gr_I^{\be}\Gr_V^{\al}A\nes0$, we say that spectral numbers $\al,\be$ are {\it shifted\one} by $\ga\defs\al\mi\be$ under the action of $[f]$, and we call $\ga$ the {\it shifted degree.} Theorem\,\,\ref{T1} shows that Tjurina spectrum is more worth studying than originally thought.
\sk
Theorem\,\,\ref{T1} implies the following.

\begin{icor} \label{C1}
There is a canonical bijection
\begin{equation*}
\phi\col\h{$\bigsqcup$}_{\ga\ges1}C'_{f,\ga}\simto C_f
\end{equation*}
with $C'_{f,\ga}\eq\bl\{1,\dots,n'_{\ga}\br\}\,\,(n'_{\ga}\ins\N)$ for $\ga\ins\Q_{\ges 1}\,\,($if $C'_{f,\ga}\nes\emptyset)$ such that
\begin{equation*}
\aligned&\phi_{\ga}(i')\slt\phi_{\ga}(j')\q\h{if}\,\,\,\,i'\slt j',\\&\phi_{\be}(i')\slt\phi_{\ga}(j')\q\h{if}\,\,\,\al_{\phi_{\be}(i')}\eq\al_{\phi_{\ga}(j')},\,\,\be\slt\ga,\\ &\al_{\phi_{\ga}(i')}\pl\al_{\phi_{\ga}(j')}\eq\ga\pl n\q\h{if}\,\,\,\,i'{+}j'\eq n'_{\ga}{+}1,\endaligned
\end{equation*}
where $\phi_{\ga}$ is the restriction of $\phi$ to $C'_{f,\ga}$ with $C_{f,\ga}\defs{\rm Im}\,\phi_{\ga}$ so that $C_f\eq\bigsqcup_{\ga\ges 1}C_{f,\ga}$.
\end{icor}

We call $C_{f,\ga}$ the shifted degree $\ga$ part of $C_f$.

\begin{irem}[{\it Shifted degree one part of missing spectral numbers\one\rm}] \label{R2}
The missing spectral numbers $\al_i$ for $i\ins C_{f,1}$ coincide with the spectral numbers of the image of $N\defs\log T_u$ in the vanishing cohomology $H^{n-1}(F_{\!f},\C)$, which can be defined by generalizing \eqref{3.3} below, where $T_u$ is the unipotent part of the Jordan decomposition $T\eq T_uT_s$ of the monodromy $T$. So the $\al_i$ for $i\ins C_{f,1}$ remain invariant under a $\mu$-constant deformation of $f$. Note that the graded morphism
\begin{equation*}
\Gr_V[f]:\Gr_V^{\ssb}A\tos\Gr_V^{\ssb+1}\!A
\end{equation*}
can be identified with the action of $-\tfrac{1}{2\pi i}\Gr_FN$ on the $F$-graded quotients of the vanishing cohomology $\Gr_F^{\ssb}H^{n-1}(F_{\!f},\C)$, see \cite{Va0}, \cite{Va1}, \cite{SS}, and also the remark after \eqref{3.1} below.
\end{irem}

Corollary\,\,\ref{C1} means that the missing spectral numbers have a canonical {\it graded symmetry.} Let $\ep$ be the involution of $C_f$ corresponding under $\phi$ to the involution of $C'_{f,\ga}$ defined by $i'\mapsto n'_{\ga}\pl1\mi i'$. Then Corollary\,\,\ref{C1} asserts that
\begin{equation} \label{2}
\al_i\pl\al_{\ep(i)}\eq\ga\pl n \q\h{if}\,\,\,i\ins C_{f,\ga}.
\end{equation}
We say that $\al_{\ep(i)}$ is a {\it graded symmetry partner\one} of a missing spectral number $\al_i$. Here the number $\ep(i)$ is not necessarily unique if the multiplicity is not 1. There are many examples of missing spectral numbers $\al_i,\al_j$ such that $\al_i\eq\al_j$ but their shifted degrees $\ga$ differ, see Section\,\,\ref{S9} below.
Note that the symmetry for the shifted degree 1 part is quite trivial using {\it Hodge symmetry\one} with the {\it hard Lefschetz property\one} for $N$, but this does not seem to suggest easily the {\it higher graded symmetry\one} where no standard Hodge theory is applicable.
Note also that the graded hard Lefschetz property under the action of $N$ never holds in general.
\sk
The missing spectral numbers do {\it not\one} have a symmetry even if $f\eq x^a\pl y^b\pl x^py^q$ $\bl(\tfrac{p}{a}\pl\tfrac{q}{b}\sgt1\br)$ unless $2p\gess a{-}1$, $2q\gess b{-}1$; for instance in the case $(a,b,p,q)\eq(7,6,5,2)$; here we have
\begin{equation*}
\mu{-}\tau\eq4\ne(a{-}1{-}p)(b{-}1{-}q)\eq3,
\end{equation*}
and the missing spectral numbers are $\tfrac{57}{42}, \tfrac{64}{42}, \tfrac{65}{42}, \tfrac{71}{42}$. They sometimes have it as a consequence of Proposition\,\,\ref{P4.2} below if the modality is very small.

\begin{irem}[{\it Quasi-weight filtration on $\Gr_V^{\ssb}{\rm Im}\,[f]$\rm}] \label{R3}
We have the quasi-weight filtration $\Wt$ indexed by $\Q$ on the graded quotients $\Gr_V^{\ssb}{\rm Im}\,[f]$ of the induced filtration $V$ on the image of the multiplication by $[f]$ in $A$ such that
\begin{equation*}
\aligned\Wt_{\ga+n}\Gr_V^{\al}{\rm Im}\,[f]&\eq I^{\al-\ga}\Gr_V^{\al}A,\\ \Gr^{\Wt}_{\ga+n}\Gr_V^{\al}{\rm Im}\,[f]&\eq\Gr_I^{\al-\ga}\Gr_V^{\al}A.\endaligned
\end{equation*}
This follows from the argument in the proof of Theorem\,\,\ref{T1}, see Section\,\,\ref{S5} below. We call $\ga{+}n$ the quasi-weight. Note that the graded vector space $\Gr^{\Wt}_{\ga+n}\Gr_V^{\ssb}$ is quite close to a complex Hodge structure of weight $e(\ga{+}n)$ replacing the indices $\al$ by $e\al$, where $e$ is the order of the semisimple part $T_s$ of the Jordan decomposition $T\eq T_sT_u$ of the monodromy $T$, however the pairing in Theorem\,\,\ref{T1} is bilinear (not {\it sesquilinear\one}) and the {\it complex conjugation\one} is missing. It seems rather interesting whether this structure can be extended to some further one.
Note that $\Gr^{\al}A$ is isomorphic to the complex conjugate of $\Gr^{n-\al}A$ using \eqref{3.1}--\eqref{3.2} and the real structure given by $H^{n-1}(F_{\!f},\R)$; however this seems useful only in the case $\ga\eq1$ (where the action of $\Gr_V[f]$ can be identified with $-\tfrac{1}{2\pi i}N$ as explained in Remark\,\,\ref{R2}), especially when $N^2\eq0$. It seems quite difficult to analyze precisely the {\it non-strict\one} part of the action of $[f]$ on the filtered ring $(A,V)$ even in the semi-homogeneous (that is, ordinary $m$-ple point) case, where there may be several shifts of degrees, for instance if $f\eq x^6\pl y^6\pl z^6\pl x^4y^4\pl x^4yz^4\pl x^2y^4z^4$ using the code in Section~\ref{S9}. (It is easy to predict the shifts of degrees $\ga$ for this kind of relatively simple examples, but this seems very difficult in general.)
\end{irem}

Corollary\,\,\ref{C1} implies for instance the following.

\begin{icor} \label{C2}
We have the inequality
\begin{equation*}
\#\{i\ins C_f\mid\al_i\slt\tfrac{n+1}{2}\}\les\tfrac{\mu-\tau}{2}.
\end{equation*}
Here the strict inequality $<$ on the left-hand side can be replaced with $\les$ if $N\eq0$ on the $\la$-eigenspace of the monodromy on the vanishing cohomology $H^{n-1}(F_{\!f},\C)_{\la}$ with $\la\defs(-1)^{n+1}$.
\end{icor}

This follows from \eqref{2}, since $\ga\ges1$ (and $|C_f|\eq\mu{-}\tau$). To show the last assertion of Corollary\,\,\ref{C2} we use the last part of Remark\,\,\ref{R2}.
This estimate seems asymptotically optimal, see Remark\,\,\ref{R5.2} below.
\sk
Using Theorem\,\,\ref{T1} and Corollary\,\,\ref{C1.1} below, we can deduce the following.

\begin{ithm} \label{T2}
The minimal missing spectral number $\al^C_1$ has multiplicity $1\,\,($that is, $\al^C_1\slt\al^C_2)$, and the maximal spectral number $\al_{\mu}\eq\al^C_{\mu-\tau}$ is the unique graded symmetry partner of $\al^C_1$.
\end{ithm}

There is an extension to the $e$-{\it Tjurina spectrum\one} case as follows.

\begin{ithm} \label{T3}
Theorems\,\,{\rm\ref{T1}--\ref{T2}} and Corollaries\,\,{\rm\ref{C1}--\ref{C2}} as well as Propositions\,\,{\rm\ref{P4.1}} and {\rm \ref{P4.2}} below can be extended to the case of $e$-Tjurina spectrum $\Spf^{(e)}(t)$ associated with the $e$-Tjurina algebra $A^{(e)}\defs\C\{x\}/(\dd f,f^e)$ of dimension $\tau^{(e)}$ for any positive integer $e\less n{-}1$. Here $C$, $\tau$, $N$, $n{+}1$, $\al_1{+}1$ must be replaced by $C^{(e)}$, $\tau^{(e)}$, $N^e$, $n{+}e$, $\al_1{+}e$ respectively, but the condition $\al{+}\be\eq\ga{+}n$ in Theorem\,\,{\rm\ref{T1}} does not change, although $\ga\gess1$ must be replaced by $\ga\gess e$.
\end{ithm}

Indeed, the same argument applies replacing the multiplication by $f$ with the one by $f^e$, see Section\,\,\ref{S9} below for calculations of explicit examples. Here $\Spf^{(1)}(t)\eq\Spf^{\rm Tj}(t)$ and
\begin{equation*}
\Spf(t)\mi\Spf^{(e)}(t)\eq\msum_{i\in\Cf^{(e)}}\,t^{\al_i}\eq\msum_{k=1}^{\mu-\tau^{(e)}}\,t^{\al_k^{C^{(e)}}}.
\end{equation*}
In the case $e\eq2$ with $n\eq3$, we have $f^e\ins(\dd f)$ if $f$ is not rather complicated (for instance if $\al_1\sgt\tfrac{1}{2}$ or $\deg f\less 6$).

\begin{irem}[{\it Composition of the action of $f\one$\rm}] \label{R4}
Theorem\,\,\ref{T3} gives some information about the ``composition" of the action of $f$ on $A$, which is quite nontrivial in general. Consider for instance the case $f\eq h{+}g_1{+}g_2$ with
\begin{equation*}
h\eq x^{11}{+}2y^{10}{+}3z^9,\q g_1\eq x^9y^2,\q g_2\eq x^4y^4z^3.
\end{equation*}
Here we have $\dim\Gr_V^{\al}A\eq1$ or 0 for any $\al\in\Q$, and a non-zero element of $\Gr_{I^{(2)}}^{\al_1}\Gr_V^{\al_1+\ga}A$ with $\al_1\eq\tfrac{299}{990}$, $\ga\eq\tfrac{2172}{990}$ is not a {\it composition\one} of elements of $\Gr_I^{\al_1}\Gr_V^{\al_1+\be}A$ and $\Gr_I^{\al_1+\be}\Gr_V^{\al_1+\ga}A$ for some $\be$. Indeed, the latter vanishes, although the former does not, for $\be\eq\tfrac{1008}{990}\eq\tfrac{9}{11}{+}\tfrac{2}{10}$, where $2\bl(\tfrac{4}{11}{+}\tfrac{4}{10}{+}\tfrac{3}{9}\br)\eq\tfrac{2172}{990}$. This can be seen by using the code in Section\,\,\ref{S9} below. It is not quite simple, since $g_1^ag_2^b\,{\notin}\,(\dd f)$ for $a,b\ins\N$ with $a{+}b\eq2$, although $g_1^2,g_1g_2\ins(\dd h)$. It seems that
\begin{equation*}
[g_1^2],[g_1g_2]\ins V^{\al}\!A\stm V^{>\al}\!A\q\h{if}\,\,\,\al\eq\tfrac{2561}{990}\eq\tfrac{299}{990}\pl\tfrac{2262}{990}.
\end{equation*}
This can be verified by computing {\smaller\sf\verb@vdim(std(I));@} for {\smaller\sf\verb@I=M[p];@} and for {\smaller\sf\verb@I=(M[p],g1^2);@} or {\smaller\sf\verb@I=(M[p],g1*g2);@} after running the first two parts of the code in Section\,\,\ref{S9} below for the above $f$, where $p=2561$ or $2562$. We have similarly $[g_2^2]\ins V^{\al}\!A\stm V^{>\al}\!A$ if $\al\eq\tfrac{2471}{990}\eq\tfrac{299}{990}\pl\tfrac{2172}{990}$.
\end{irem}

Related to an estimate of the {\it Brian\c{c}on-Skoda exponent\one} $e^{\rm BS}(f)$ in \cite[Theorem 1.3]{Va1} (whose proof is simplified very much and is extended to the non-isolated singularity case in \cite{JKSY1} avoiding the asymptotic Hodge filtration and using the $V$-filtration of Kashiwara and Malgrange indexed by $\Q$, which was originally indexed by $\Z$), set
\begin{equation*}
\aligned\si_{\!f}&\defs\min\bl\{\ga\ins\Q\,\,\big|\,\,\Gr_I^{\al-\ga}\Gr_V^{\al}A\ne0\,\,(\exists\,\al\ins\Q)\br\},\\ \si'_{\!f}&\defs\max\bl\{\ga\ins\Q\,\,\big|\,\,
[f]V^{\be}\!A\subset V^{\be+\ga}\!A\,\,(\forall\,\be\ins\Q)\br\}.\endaligned
\end{equation*}
We have the following.

\begin{ithm} \label{T4}
There is the equality
\begin{equation} \label{3}
\si_{\!f}=\si'_{\!f},
\end{equation}
and we have the inequality
\begin{equation*}
e^{\rm BS}(f)\les\Bigl[\frac{n{-}2\al_1}{\si_{\!f}}\Bigr]\pl1,
\end{equation*}
or equivalently, $f^k\ins(\dd f)$ if $\,k\one\si_{\!f}\sgt n\mi2\al_1$.
\end{ithm}

The last assertion follows from the first using the {\it equality\one} \eqref{1} together with $\al_{\mu}\eq n\mi\al_1$ in \eqref{3.4}, see also Proposition\,\,\ref{P1.1} below.
(Note that \cite[Theorem 2]{JKSY1} in the {\it isolated\one} singularity case is a consequence of the coincidence of the two quotient $V$-filtrations on the Jacobian ring.)
This theorem improves \cite[Theorem 1.3]{Va1} or \cite[Theorem 1]{JKSY1} for {\it isolated\one} singularities {\it only\one} when the monodromy $T$ is {\it semisimple,} that is, $N\eq0$, or equivalently $\si'_{\!f}\sgt1$. (The proof of \cite[Theorem 1.3]{Va1} is reduced to a quite nontrivial assertion, that is, Theorem 9.1, whose proof seems very complicated due to the difference between the Steenbrink and asymptotic Hodge filtrations; for instance, it does not seem entirely trivial that the leading term $s_{\max}[\om'](t)$ of the asymptotic expansion of $[\om']\ins\Hsf$ belongs to the subspace $W_l$, where some more explanation seems desirable.) Set
\begin{equation*}
\si''_{\!f}\defs\max\{\ga\ins\Q\,|\,[f]\ins V^{\al_1+\ga}\!A\}.
\end{equation*}
We have the following.

\begin{iprop} \label{P1}
In the semi-weighted-homogeneous case, there are equalities
\begin{equation} \label{4}
\si_{\!f}\eq\si'_{\!f}\eq\si''_{\!f},
\end{equation}
and we have $\si_{\!f}\gess\min\{\be\sgt1\mid f_{\be}\nes0\}$ if $f\eq\msum_{\be\ges1}\,f_{\be}$ in the notation of Remark\,\,{\rm\ref{R3.5}} below. Moreover the dimension of $\Gr_I^{\al-\ga}\Gr_V^{\al}A$ does not change by adding $g\ins V_{\mhs N}^{>\ga+\al_1}\C\{x\}$ to $f$, and this assertion holds also in the convenient Newton non-degenerate case for $g\ins V_{\mhs N}^{\prime\one >\ga}\C\{x\}$ in the notation of Remark\,\,{\rm\ref{R3.4}} below.
\end{iprop}

The last equality of \eqref{4} does not necessarily hold in the Newton non-degenerate case, for instance if the intersection of the line spanned by $(1,\dots,1)$ with the Newton boundary is contained in the interior of a maximal-dimensional compact face of the Newton polytope and the monodromy is {\it not\one} semisimple, where $\si''_{\!f}\sgt\si_{\!f}\eq1$, see Remark\,\,\ref{R4.1} below.
\sk
We thank the referees for careful reading and useful comments on this paper.

\tableofcontents
\numberwithin{equation}{section}

\section{Brieskorn lattices and Gauss-Manin systems} \label{S1}
Let $f\ins\C\{x\}$ be a convergent power series having an isolated singularity with $x\eq(x_1,\dots,x_n)$ the natural coordinate system of $(\C^n,0)$ and $f(0)\eq0$. The {\it Brieskorn lattice\one} (see \cite{Br}) is defined by
\begin{equation*}
\Hsf\defs\Om_{\C^n,0}^n/\ddd f{\wedge}\one\ddd\Om_{\C^n,0}^{n-2}.
\end{equation*}
This is a free module of rank $\mu$ over $\C\{t\}$ and also over $\C\{\!\{\dti\}\!\}$, where $\mu$ is the Milnor number and
\begin{equation*}
\Rc:=\C\{\!\{\dti\}\!\}=\bl\{\msum_{j\in\N}\,c_i\dt^{-j}\ins\C[[\dti]]\,\,\big|\,\,\msum_{j\in\N}\,c_jt^j\!/j!\ins\C\{t\}\br\},
\end{equation*}
see for instance \cite{bl}, \cite{mic}. The actions of $t$, $\dti$ on $\Hsf$ are defined by
\begin{equation} \label{1.1}
\aligned&t[\om]\eq[f\om],\q\dti[\om]\eq[\ddd f{\wedge}\one\eta],\\&\q\h{where}\,\,\,\om\ins\Om_{\C^n,0}^n,\,\eta\ins\Om_{\C^n,0}^{n-1}\,\,\,\h{satisfying}\,\,\,\ddd\eta\eq\om.\endaligned1.
\end{equation}
\sk
Set $\K\defs \Rc[\dt]$ (the localization of $\Rc$ by $\dti$). The {\it Gauss-Manin system\one} $\Gf$ is isomorphic to the localization of the Brieskorn lattice $\Hsf$ by the action of $\dti$ so that
\begin{equation*}
\Gf=\Hsf[\dt].
\end{equation*}
This is a free $\K$-module of rank $\mu$, and has the $V$-filtration indexed by $\Q$ such that the $V^{\al}\Gf$ are free $\Rc$-modules of rank $\mu$ and satisfying
\begin{equation} \label{1.2}
\aligned t\one V^{\al}\Gf&\sst V^{\al+1}\Gf,\\ \dt\one V^{\al}\Gf&\eq V^{\al-1}\Gf,\\(\dt t\mi\al)^n\one &\Gr_V^{\al}\Gf\eq0\q\q(\forall\one\al\ins\Q).\endaligned
\end{equation}
\sk
There is the Hodge filtration $F$ on $\Gf$ defined by
\begin{equation} \label{1.3}
F\!_p\one\Gf\defs\dt^{\one p+n-1}\Hsf\q(\forall\,p\ins\Z),
\end{equation}
see \cite[2.7]{bl} (where the number of variables is $n{+}1$ instead of $n$ as in our paper). This is shifted by $1$ compared with the right $\D$-module case as in \cite{mhp}. (Note that one gets the shift by $1$ taking $\Gr_V^{\al}$ for $\al\ins(0,1]$ in the right $\D$-module case.) Set
\begin{equation*}
\Om_f\defs\Om_{\C^n,0}^n/\ddd f{\wedge}\Om_{\C^n,0}^{n-1}.
\end{equation*}
By definition we have the isomorphisms
\begin{equation} \label{1.4}
\Gr^F_{1-n}\Gf=\Hsf/\dti\mhs\Hsf=\Om_f\cong A\,\bl({:=}\,\C\{x\}/(\dd f)\br),
\end{equation}
where the last isomorphism depends on the generator $\ddd x\defs\ddd x_1{\wedge}\cdots{\wedge}\ddd x_n$ of $\Om_{\C^n,0}^n$. (In this paper equalities and canonical isomorphisms are denoted by $=$, and non-canonical ones by $\cong$ basically.)
\sk
From \eqref{1.1}--\eqref{1.3} we can deduce the following important assertion (which immediately implies an improvement of the Brian\c{c}on-Skoda theorem \cite{BrSk} in the {\it isolated\one} singularity case simplifying the proof of \cite[Theorem 1.3]{Va1} as explained after Theorem\,\,\ref{T4}).

\begin{prop} \label{P1.1}
For any $\al\ins\Q$, we have the inclusion
\begin{equation*}
f\one V^{\al}\bl(\C\{x\}/(\dd f)\br)\sst V^{\al+1}\bl(\C\{x\}/(\dd f)\br).
\end{equation*}
\end{prop}

We have also the following.

\begin{prop} \label{P1.2}
The $V$-filtration on the Jacobian ring $A\eq\C\{x\}/(\dd f)$ is a filtration by ideals, and its graded pieces are annihilated by the maximal ideal induced by that of $\C\{x\}$.
\end{prop}

This is proved using the coincidence of the $V$-filtration on the Jacobian ring with the quotient filtration of the microlocal $V$-filtration on $\C\{x\}$, see \cite[Rem.\,3.11]{pmb}, \cite[\S4.11]{DS2}, \cite{JKSY}.
\sk

Employing the {\it symmetry of spectral numbers\one} (see \eqref{3.4} below) and the {\it equality\one}~\eqref{1}, we can easily deduce the following.

\begin{cor}\label{C1.1}
The filtered $\C$-algebra $A$ is a filtered Artinian local ring with maximal ideal $\mm\sst A$ induced by that of $\C\{x\}$, and we have
\begin{align}\label{1.5}
V^{\al_1}A\eq A,\q V^{>\al_1}A\eq\mm,&\q V^{>\al_{\mu}}A\eq 0,\q\Gr_V^{\al_1}A\cong\Gr_V^{\al_{\mu}}A\cong\C,\\
[f]A\subset V^{\al_1+1}A,&\q[f]\one V^{>\al_1}A\eq[f]\mm\subset V^{>\al^C_1}A,\label{1.6}
\end{align}
where $[f]\ins A$ denotes the class of $f$.
\end{cor}

(Since the argument is quite straightforward, the details are left to the reader.)

\section{Self-duality of bifiltered Gauss-Manin systems} \label{S2}
We have the following.

\begin{thm}[{\cite[2.7]{bl}}] \label{T1.2}
There is a canonical pairing
\begin{equation*}
\Sb:\Gf\otimes_{\C}\Gf\to\K\,\bl({=}\,\C\{\!\{\dti\}\!\}[\dt]\br),
\end{equation*}
satisfying
\begin{equation} \label{2.1}
\aligned P\one\Sb(u,v)&\eq\Sb(Pu,v)\eq\Sb(u,P^*v),\\t\one\Sb(u,v)&\eq\Sb(tu,v)\mi\Sb(u,tv),\\\Sb(F_p\Gf,&F_q\Gf)\subset F_{p+q+n-2}\one\K,\endaligned
\end{equation}
for any $u,v\ins\Gf$, $P\ins\K$, $p,q\ins\Z$ with $(\dt^j)^*\defs(-1)^j\dt^j$, and the induced pairing
\begin{equation} \label{2.2}
\Gr^F\Sb:\Gr^F_p\Gf\otimes_{\C}\Gr^F_q\Gf\to\C\one\dt^{\one p+q+n-2}\,\bl({=}\,\Gr^F_{p+q+n-2}\one\K\br)
\end{equation}
is non-degenerate, where the filtration $F$ on $\K$ is by the order of $\dt$. Moreover the pairing $\Sb$ is strictly compatible with the $V$-filtration so that
\begin{equation} \label{2.3}
\Sb(V^{\al}\Gf,V^{\be}\Gf)\subset\Rc\dt^{-\lceil\al+\be\rceil}\q(\forall\,\al,\be\ins\Q),
\end{equation}
with $\Rc\eq\C\{\!\{\dti\}\!\}$, and the induced pairing
\begin{equation} \label{2.4}
\Gr^F\Gr_V\Sb:\Gr^F_p\Gr_V^{\al}\Gf\otimes_{\C}\Gr^F_q\Gr_V^{\be}\Gf\to\C\one\dt^{-k}\,\bl({=}\,\Gr^F_{-k}\one\K\br)
\end{equation}
for $-p\mi q\mi n\pl 2\eq\al\pl\be\eq k$ is non-degenerate.
\end{thm}

\begin{rem}[{\it Duality of the filtrations\one\rm}] \label{R2.1}
As a corollary of Theorem\,\,\ref{T1.2}, the filtration $V^{\ssb}$ on $A\eq\C\{x\}/(\dd f)$ is the dual of $V^{\ssb}[n]$ with respect to $\Gr^F\Sb$, that is,
\begin{equation} \label{2.5}
V^{>n-\al}\!A\eq(V^{\al}\!A)^{\perp}\defs\bl\{v\ins A\mid\Gr^F\Sb(v,v')\eq0\,\,(\forall\,v'\ins V^{\al}\!A)\br\}.
\end{equation}
We then see that
\begin{equation} \label{2.6}
(I^{\al})^{\perp}\eq\bl\{v\ins A\mid\Gr^F\Sb([f]v,v')\eq0\,\,(\forall\,v'\ins V^{\al}\!A)\br\}\eq(V^{>n-\al}\!A\col[f]),
\end{equation}
where $(B\col[f])\eq\{v\ins A\mid [f]v\ins B\}$ for $B\sst A$ in general. Indeed, we have
\begin{equation} \label{2.7}
\Gr^F\Sb([f]v,v')=\Gr^F\Sb(v,[f]v')\q(\forall\,v,v'\ins A),
\end{equation}
by \eqref{2.1} and \eqref{1.1}. Note that the action of $t$ on $\K$ is defined by $t\dd_t^j\eq{-}j\dd_t^{j-1}$ for $j\ins\Z$. (Indeed, $\K$ is isomorphic to $\E/\E t$ as $\E$-module with $\E$ the ring of microdifferential operators, so it has a generator {\it annihilated by\one} $t$.)
\end{rem}

\begin{rem}[{\it Some technical difficulty\one\rm}] \label{R2.2}
For the last part of Theorem\,\,\ref{T1.2}, we use the assertion that the canonical pairing induces a polarization of mixed Hodge structures via the isomorphisms \eqref{3.1}. It does not seem quite clear whether the non-degeneracy of \eqref{2.4} follows from the remaining part, although it does forgetting $F$. Indeed, by \eqref{2.1} we get
\begin{equation} \label{2.8}
\dt t\one\Sb(u,v)\eq\Sb(\dt tu,v)\pl\Sb(u,\dt tv).
\end{equation}
In the case the monodromy is semisimple, this implies that
\begin{equation} \label{2.9}
\Sb(\Gf^{\al},\Gf^{\be})\sst\C\one\dt^{-\al-\be}\,\,\,\h{if}\,\,\,\al{+}\be\ins\Z,\,\,\,\h{and 0 otherwise,}
\end{equation}
where $\Gf^{\al}\defs{\rm Ker}(\dd_tt\mi\al)^k\sst\Gf$ ($k\,{\gg}\,0$).
In general we can use a filtration of $\Gf$ such that the graded pieces are simple $\K\langle t\rangle$-modules together with the assertion that $\Gf$ is a regular holonomic $\K\langle t\rangle$-module, that is, $\Gf$ is generated over $\K$ by the $\Gf^{\al}$.
\end{rem}

\begin{rem}[{\it Grothendieck residue pairing\one\rm}] \label{R2.3}
In the case $p\eq q\eq 1{-}n$, it is known (see \cite{bl}) that the induced self-pairing $\Gr^F\Sb$ on $\Gr^F_{1-n}\Gf\eq\Om_f$ can be identified with the {\it Grothendieck residue pairing\one}:
\begin{equation} \label{2.10}
{\rm Res}_0\Bigl[\begin{array}{c}gh\ddd x\\f_1\cdots f_n\end{array}\Bigr]=\frac{1}{(2\pi i)^n}\int_{|f_1|\eq\ep_1,\,\dots\,,|f_n|=\ep_n}\frac{gh\ddd x}{f_1\cdots f_n}
\end{equation}
for $g,h\ins\C\{x\}$, where $f_i$ denotes $\dd_{x_i}f$ in this Remark, and $0\slt\ep_i\,{\ll}\,1$, see \cite[p.\,195]{Ha}, \cite[p.\,659]{GH}.
\end{rem}

\begin{rem}[{\it A conjecture of K. Saito\one\rm}] \label{R2.4}
It has been conjectured by K.\,Saito that the difference between the maximal and minimal exponents of a hypersurface isolated singularity is at most 1 if and only if the singularity is rational double or simple elliptic or cusp, see \cite{SWZ}. Here the exponents are defined by using ``good sections", and do not necessarily coincide with the spectral numbers defined by Steenbrink \cite{St} using Hodge theory, see for instance \cite{bl}, \cite{bl2}, \cite{JKSY3}. It seems desirable to eliminate bad examples by adding some appropriate condition, although using $V$-filtration as in this paper does not seem quite adequate for the theory.
\end{rem}

\section{Steenbrink spectrum} \label{S3}
In the notation of Section\,\,\ref{S1}, there are canonical isomorphisms
\begin{equation} \label{3.1}
H^{n-1}(F_{\!f},\C)_{\la}=\Gr_V^{\al}\Gf\q(\forall\,\al\ins(0,1],\,\la\eq e^{-2\pi i\al}),
\end{equation}
with coordinate $t$ of a disk $\De$ fixed. Here $F_{\!f}$ denotes the Milnor fiber of $f$, and $H^{n-1}(F_{\!f},\C)_{\la}$ means the $\la$-eigenspace for the semisimple part of the monodromy (which is the inverse of the Milnor monodromy, see \cite{DS1}). Moreover the action of $\dd_tt\mi\al$ on the right-hand side corresponds to $-\tfrac{1}{2\pi i}N$ on the left-hand side.
These follow from the theory of Deligne extensions of local systems. Indeed, a multivalued horizontal section $u$ of a $\C$-local system $L$ on a punctured disk $\De^*$ corresponds to
\begin{equation*}
\widetilde{u}\defs t^{\al-1}e^{-(2\pi i)^{-1}N\log t}u\q\h{in}\,\,\,\,j_*(\OO_{\De^*}{\otimes}_{\C}L),
\end{equation*}
with $j\col\De^*\into\De$ the inclusion, if $u$ is annihilated by $T_s{-}e^{-2\pi i\al}$ with $T_s$ the semisimple part of the Jordan decomposition $T\eq T_sT_u$ of the monodromy $T$ (and $N\defs\log T_u$).
\sk
Let $F$ be the Hodge filtration on the vanishing cohomology $H^{n-1}(F_{\!f},\C)$ (see \cite{St}), which is compatible with the decomposition by the eigenspaces of the semisimple part $T_s$ of the monodromy. This Hodge filtration coincides with the Hodge filtration on the vanishing cycle Hodge module, which is denoted by $\varphi_f\Q_{h,X}[n{-}1]$ in this paper, see \cite{mhm}. One can verify this by employing a compactification of $f$ as in \cite{SS} (see also \cite{gm1}), since the vanishing cycle functor commutes with the cohomological direct image functor under a projective morphism. Using this, one can show for instance the invariance of the spectrum by a $\mu$-constant deformation of $f$, which was obtained in \cite{Va2}.
\sk
It is quite well known that the isomorphism \eqref{3.1} induces the isomorphisms
\begin{equation} \label{3.2}
F^p\!H^{n-1}(F_{\!f},\C)_{\la}=F\!_{-p}\Gr_V^{\al}\Gf\q(\forall\,p\ins\Z,\,\al\ins(0,1],\,\la\eq e^{-2\pi i\al}),
\end{equation}
see \cite{SS}, \cite{Va1} (and also \cite[\S3.4]{gm1}).
\sk
In the notation of the introduction, the Steenbrink spectrum of $f$ is usually defined by
\begin{equation} \label{3.3}
|\If^{\al}|=\dim_{\C}\Gr_F^pH^{n-1}(F_{\!f},\C)_{\la}\q(p\defs[n{-}\al],\,\la\defs e^{-2\pi i\al}).
\end{equation}
This is equivalent to the definition \eqref{1} in the introduction by \eqref{3.2}.
\sk
We can prove the {\it symmetry\one} of spectral numbers using Hodge theory (see \cite{St}):
\begin{equation} \label{3.4}
\al_i\pl\al_j\eq n\q\h{if}\,\,\,i{+}j\eq\mu{+}1.
\end{equation}
Here one needs in an essential way the assertion that the weight filtration is given by the shifted monodromy filtration (where the shift depends on the monodromy eigenvalues).

\begin{rem}[{\it Original formulation\one\rm}] \label{R3.1}
In the original definition of spectrum (see \cite[\S5.3]{St}) the formulation is slightly different, since $q$ is used for the integer part (although it is not shifted by $-1$ as in \cite{SS}). This is closely related to a confusion about the relation between the monodromy and the Milnor monodromy, which is explained in \cite[\S2.1]{DS1}.
\end{rem}

\begin{rem}[{\it Another normalization\one\rm}] \label{R3.2}
In some papers, the spectrum is shifted by $-1$, for instance in \cite{SS}, \cite{Sing}. This is closely related to the asymptotic expansions of period integrals \cite{Va1}. It is better to employ the unshifted definition in case one considers the relation to the multiplier ideals or to the roots of Bernstein-Sato polynomials.
\end{rem}

\begin{rem}[{\it Formula in the weighted homogeneous case\one\rm}] \label{R3.3}
There is a well-known explicit formula for the Steenbrink spectrum in the weighted homogeneous case:
\begin{equation} \label{3.5}
\Spf(t)\eq\mprod_{i=1}^n\,(t^{w_i}\mi t)/(1\mi t^{w_i}),
\end{equation}
where the $w_i$ are the weights of variables $x_i$ such that $\msum_{i=1}^n\,w_ix_i\dd_{x_i}f\eq f$. This was found by Steenbrink \cite{St} without giving any proof. Its proof is not quite trivial unless all the weights are inverses of integers (where the assertion can be reduced to the Brieskorn-Pham polynomial case), see for instance \cite{JKY}, \cite{JKSY}.)
\end{rem}

\begin{rem}[{\it Newton non-degenerate case\one\rm}] \label{R3.4}
Assume $f$ is Newton non-degenerate and also convenient (that is, for any $i\ins[1,n]$, the coefficient of $x_i^{a_i}$ in $f$ does not vanish for some $a_i\ins\N$). We have the {\it Newton filtrations\one} $V'_{\mhs N}$, $V_{\mhs N}$ on $\C\{x\}$ consisting of ideals generated by monomials and such that
\begin{equation} \label{3.6}
\aligned x^{\nu}\ins V_{\mhs N}^{\prime\one\al}\C\{x\}&\iff\nu\in\al\one\Ga_{\!+}(f)\q(\forall\,\al\ins\Q_{>0}),\\ x^{\nu}\ins V_{\mhs N}^{\al}\C\{x\}&\iff x^{\nu+{\bf 1}}\ins V_{\mhs N}^{\prime\one\al}\C\{x\}\q(\forall\,\al\ins\Q_{>0}),\endaligned
\end{equation}
where ${\bf 1}\defs(1,\dots,1)\ins\R^n$, see \cite{Ko}, \cite{exp}, \cite{JKSY2}. These filtrations are exhaustive since $f$ is convenient. Note that $(\C\{x\},V'_{\mhs N})$ is a filtered {\it ring,} and $(\C\{x\},V_{\mhs N})$ is a filtered {\it module\one} over it. As a corollary of the proof of the Steenbrink formula in \cite{JKSY2}, we get that
\begin{equation} \label{3.7}
\aligned&\h{The quotient filtration of $V_{\mhs N}$ on $\C\{x\}/(\dd f)$ coincides with the quotient}\\
&\h{filtration of the $V$-filtration on the Brieskorn lattice $\Hsf$ via \eqref{1.4}.}\endaligned
\end{equation}
Here \cite{exp} does not seem quite sufficient. Indeed, it proves only an assertion on the $V$-filtration on the Gauss-Manin system, and strictly speaking, some additional argument seems to be required actually, see \cite[Remark 2.1f]{JKSY2} and also \cite[Remark 2.2d]{exa}.
\end{rem}

\begin{rem}[{\it Semi-weighted-homogeneous case\one\rm}] \label{R3.5}
A polynomial $f$ is called {\it semi-weighted-homogeneous\one} if $f\eq\msum_{\be\ges1}f_{\be}$ with $f_{\be}$ weighted homogeneous of degree $\be$, that is, $\xi_{\bf w}f_{\be}\eq\be f_{\be}$ with $\xi_{\bf w}\defs\msum_iw_ix_i\dd_{x_i}$ ($\be\gess 1$), and $f_1$ has an isolated singularity at 0. In this case we have by \cite{Va2}
\begin{equation} \label{3.8}
\Spf(t)\eq\Sp_{\mhs f_1}(t).
\end{equation}
Moreover \eqref{3.7} holds also in this case, that is, the $V$-filtration $V^{\al}$ on $\C\{x\}/(\dd f)$ is generated by the monomials $x^{\nu}$ with
\begin{equation} \label{3.9}
\msum_iw_i(\nu_i\pl1)\eq\msum_iw_i\nu_i\pl\al_1\ges\al,
\end{equation}
One can easily verify this in the case $f\eq f_1$ using the Euler field $\xi_{\bf w}$, but the argument in the general semi-weighted-homogeneous case does not seem quite trivial. It is easy to prove an inclusion, and the coincidence of the dimensions of the graded pieces can be verified by comparing the graded quotients of the Koszul complexes, see also \cite[Remark 1.2d]{len}, \cite[Remark 2.2d]{exa} (and Remark\,\,\ref{R3.4} in the case $1/w_i\ins\Z$ for any $i$).
\end{rem}

\begin{rem}[{\it Saturated Hodge filtration\one\rm}] \label{R3.6}
We can define the {\it saturated Hodge filtration\one} $\Ft$ on the vanishing cohomology $H^{n-1}(F\!_f,\C)_{\la}$ using the isomorphism \eqref{3.1} and replacing the Brieskorn lattice $\Hsf\sst\Gf$ with its {\it saturation\one} $\widetilde{H}''_{\!f}\defs\msum_{j\ges0}\,(\dd_tt)^j\Hsf\sst\Gf$, see \cite[\S8]{Va1} (and also \cite{sem}). Recall that the reduced Bernstein-Sato polynomial $\widetilde{b}_f(s)\defs b_f(s)/(s{+}1)$ is equal to the minimal polynomial of the action of $-\dd_tt$ on $\widetilde{H}''_{\!f}/t\widetilde{H}''_{\!f}$, see \cite{Ma}. We define the $b$-exponents by replacing the Hodge filtration $F$ with the {\it saturated\one} one $\Ft$, where the minimal $b$-exponent coincides with the minimal spectral number $\al_1$ by definition. The ``shifts between the exponents" can be seen from the bigraded pieces $\Gr_F^p\Gr_{\Ft}^qH^{n-1}(F\!_f,\C)_{\la}$. These imply immediately that any spectral number which is strictly smaller than $\al_1{+}1$ is a root of $b_f(s)$ up to sign.
\sk
The $b$-exponents were used in a conjecture of Yano for the irreducible plane curve case and the related invariants, which were defined by taking the pullback of the highest differential form to the embedded resolution, were calculated in the handwritten version of \cite{irr}, 1982.
\end{rem}

\begin{rem}[{\it Thom-Sebastiani type situation\one\rm}] \label{R3.7}
Assume $X\eq X'{\times}X''$, $f\eq{\rm pr}'{}^*f'\pl{\rm pr}''{}^*f''$ with ${\rm pr}'\col X\tos X'$ the projection (similarly for ${\rm pr}''$). Let $\mu',\mu''$ and $\tau',\tau''$ be the Milnor and Tjurina numbers of $f',f''$ respectively. It is rather difficult to determine $\tau$ from $\tau',\tau''$, see for instance \cite{Al}. This is the same for the Tjurina spectrum even though we have the compatibility of the following Thom-Sebastiani type isomorphism with the $V$-filtration (see \cite{SS}, \cite{Va1}):
\begin{equation} \label{3.10}
\C\{x\}/(\dd f)=\bl(\C\{x'\}/(\dd f')\br){\otimes}\bl(\C\{x''\}/(\dd f'')\br).
\end{equation}
This filtered isomorphism however implies the strict compatibility with the $V$-filtration of the canonical surjection
\begin{equation} \label{3.11}
\C\{x\}/(\dd f,f)\onto\bl(\C\{x'\}/(\dd f',f')\br){\otimes}\bl(\C\{x''\}/(\dd f'',f'')\br).
\end{equation}
So $\SpT_{\mhs f'}(t)\one {\cdot}\one\SpT_{\mhs f''}(t)$ is a subspectrum of $\SpTf(t)$, that is, the difference has positive coefficients. We can verify that the monomial whose exponent is the spectral number
\begin{equation*}
\max\bl(\al'_1\pl\al_1^{\prime\prime\one C},\,\al''_1\pl\al_1^{\prime\one C}\br)
\end{equation*}
always belongs to this difference and the minimum of the above two numbers is a missing spectral number of $f$ looking at the maximal number $\al$ such that the image of $f$ in the right hand side of \eqref{3.10} (that is, $f'{\otimes}1\pl1{\otimes}f''$) is contained in $V^{\al}$. Here $\al'_1$ and $\al_1^{\prime\one C}$ are the minimal (missing) spectral number of $f'$, and similarly for $f''$. In the simplest case with $\mu'{-}\tau'\eq\mu''{-}\tau''\eq1$, the difference consists only of this monomial.
\end{rem}

\section{Some related propositions} \label{S4}

We explain some assertions on Tjurina spectrum, which are helpful to understand the main theorems. (Some of them are moved from an earlier version of \cite{JKSY3}.) Applying the self-duality of the Jacobian ring used in the proof of Theorem\,\,\ref{T1}, we can give a proof of the following (which has been conjectured in \cite{SWZ}, and would be useful to simplify some arguments there).

\begin{prop} \label{P4.1}
In the case $\mu\nes\tau$, the maximal spectral number $\al_{\mu}$ is always a missing spectral number, that is, $\mu\ins\Cf$, hence $\al^C_{\mu-\tau}\eq\al_{\mu}\eq n\mi\al_1$.
\end{prop}

\begin{proof}
We have the endomorphism $\rho\ins{\rm End}_{\C}(\Om_f)$ defined by the action of $f$ on $\Om_f$, which is identified with $A$ using $\ddd x$. This endomorphism is self-dual for the induced perfect self-pairing of $\Om_f$ in \eqref{2.2}, see \eqref{2.7}. So we get the induced perfect pairing
\begin{equation} \label{4.1}
(\Om_f/{\rm Ker}\,\rho)\otimes_{\C}{\rm Im}\,\rho\to\C.
\end{equation}
This is strictly compatible with the $V$-filtration as a corollary of the non-degeneracy of the induced pairing \eqref{2.4}. Note that $k\eq n$ if $p\eq q\eq 1{-}n$. We then get the induced perfect pairings
\begin{equation} \label{4.2}
\Gr_V^{\al}(\Om_f/{\rm Ker}\,\rho)\otimes_{\C}\Gr_V^{n-\al}{\rm Im}\,\rho\to\C\q(\forall\,\al\ins\Q).
\end{equation}
The assertion now follows from Proposition\,\,\ref{P1.2}. Indeed, we see that ${\rm Ker}\,\rho$ is contained in the maximal ideal of $A$, since otherwise ${\rm Ker}\,\rho$ coincides with $A$ (and $\mu\nes\tau$). This implies that $\Gr_V^{\al_1}(\Om_f/{\rm Ker}\,\rho)$ does not vanish, hence $\Gr_V^{n-\al_1}{\rm Im}\,\rho\nes0$ by the perfectness of the pairing. This terminates the proof of Proposition\,\,\ref{P4.1}.
\end{proof}

We have also the following.

\begin{prop} \label{P4.2}
In the case $\mu\nes\tau$, the minimal and maximal missing spectral numbers $\al^C_1$ and $\al^C_{\mu-\tau}$ have both multiplicity one in $\Sp^C_f(t)$, and
\begin{equation} \label{4.3}
\al^C_1\gess\al_1\pl1,\q\h{that is,}\q\Cf\sst\If^{\ges\al_1+1},
\end{equation}
with $\If^{\ges\be}\defs\{i\ins\If\mid\al_i\gess\be\}$ for $\be\ins\Q$. Moreover the strict inequality holds in {\rm\eqref{4.3}} if the monodromy is semisimple or more generally if $\Gr_V^{\al_1}[\ddd x]\ins\Gr_V^{\al_1}\Hsf$ is annihilated by $\dd_tt\mi\al_1$.
\end{prop}

\begin{proof}
The first two assertions follow from Corollary\,\,\ref{C1.1}, and the last one from the last remark in Remark\,\,\ref{R4.1} just below. This finishes the proof of Proposition\,\,\ref{P4.2}.
\end{proof}

\begin{rem}[{\it Some semisimplicity argument\one\rm}] \label{R4.1}
In the convenient Newton non-degenerate case, the last hypothesis of Proposition\,\,\ref{P4.2} is satisfied if the intersection of the line spanned by ${\bf 1}\defs(1,\dots,1)\ins\R^n$ with the boundary of the Newton polytope is contained in the interior of some maximal-dimensional compact face $\si_0$. Indeed, the latter face is contained in a hyperplane defined by
\begin{equation*}
\msum_{i=1}^m\,c_i\nu_i=1\q\h{for some}\,\,\,\,c_i\ins\Q_{>0}\,\,(i\ins[1,n]),
\end{equation*}
and $\al_1\eq\msum_{i=1}^n\,c_i$ by the assumption. We can verify that
\begin{equation} \label{4.4}
f\mi\msum_{i=1}^n\,c_ix_if_i\ins V_{\mhs N}^{>\al_1+1}\C\{x\},
\end{equation}
where $f_i$ denotes $\dd_{x_i}f$ in this Remark. Indeed, setting
\begin{equation*}
\ell(\nu):=\min\{\ell_{\si}(\nu)\},
\end{equation*}
with $\ell_{\si}$ the linear function on $\R^n$ such that $\si\sst\ell_{\si}^{-1}(1)$ for maximal-dimensional compact faces $\si$, we can see that
\begin{equation} \label{4.5}
\ell(\nu{+}{\bf 1})\sgt1{+}\al_1\q\h{for any}\,\,\,\,\nu\in\Ga_{\!+}(f)\stm\si_0,
\end{equation}
studying the derivative of the piecewise linear function $\phi(s)\defs\ell(\nu{+}s{\bf 1})$ for $s\ins[0,1]$.
\sk
The assertion \eqref{4.4} implies the inequality $\si''_{\!f}\sgt1$ and also that
\begin{equation} \label{4.6}
\aligned(t\mi\al_1\dti)[\ddd x]&\in V^{>\al_1+1}\Hsf,\,\,\,\h{that is,}\\(\dd_tt\mi\al_1)[\ddd x]&=0\,\,\,\h{in}\,\,\,\Gr_V^{\al_1}\Gf.\endaligned
\end{equation}
The last assertion is compatible with a formula for spectral pairs in \cite{JKSY2}.
\sk
Note that the equivalence in \eqref{4.6} holds in the general isolated hypersurface singularity case, and the first condition implies that $[f]\ins V^{>\al_1+1}A$.
\end{rem}

We provide also a proof of the following.

\begin{prop} \label{P4.3}
The $V$-filtration on the Tjurina algebra $\C\{x\}/(\dd f,f)$ is independent of a defining function $f$ of the hypersurface, that is, the Tjurina spectrum depends only on the analytic hypersurface germ $(f^{-1}(0),0)\sst(\C^n,0)$.
\end{prop}

\begin{proof}
One can easily verify that the Tjurina ideal $(\dd f,f)\sst\C\{x\}$ is independent of a choice of defining function $f$. The $V$-filtration on $\C\{x\}/(\dd f,f)$ can be induced from the {\it microlocal\one} $V$-filtration on $\C\{x\}$, see \cite[\S4.11]{DS2}, \cite{JKSY}. Moreover the latter coincides with the Hodge ideals module $(f)$, see \cite[\S2.4]{JKSY} (where the coefficients $\al_k$ are independent of $k$ in the case $n\eq2$) and also \cite{him} for the case $\al_k\eq1$. So the assertion follows. This completes the proof of Proposition\,\,\ref{P4.3}.
\end{proof}

\begin{rem}\label{R4.2}
Proposition\,\,\ref{P4.3} was essentially shown in the proof of the assertion claiming that the Tjurina spectrum is a subspectrum of the Steenbrink spectrum and also of the {\it Hodge ideal spectrum\one} (see the introduction of \cite{JKSY} and also \cite[\S2.4]{SWZ}) although it does not seem to be stated explicitly there. The above proof of Proposition\,\,\ref{P4.3} employs the {\it analytic version\one} of a theorem of Musta\c{t}\u{a} and Popa \cite{MP} on the coincidence of the Hodge ideals with the {\it microlocal\one} $V$-filtration modulo the ideal of the hypersurface (see \cite[\S2.4]{JKSY} and also \cite{him} for the integer index case) using the coincidence of two quotient filtrations on the Jacobian ring which are induced by the microlocal $V$-filtration on the structure sheaf and by the $V$-filtration on the Brieskorn lattice (\cite[\S4.11]{DS2}, \cite{JKSY}). It does not seem trivial to demonstrate for instance the independence of the analytic Hodge ideals under the choice of algebraization {\it without using the analytic theory of Hodge ideals as in \cite{JKSY}.})
\end{rem}

\begin{rem}[{\it Spectral numbers and roots of Bernstein-Sato polynomials\one\rm}] \label{R4.3}
It is well known that any spectral number which is strictly smaller than $\al_1\pl1$ is a root of the Bernstein-Sato polynomial $b_f(s)$ up to sign. This is a corollary of \cite{Ma} and \cite{Va1}, see Remark\,\,\ref{R3.6} below. In particular, it holds also for Tjurina spectral numbers. It is not necessarily true that any Tjurina spectral number is a root of $b_f(s)$ up to sign; for instance if $f\eq x^9\pl y^7\pl x^4y^4$, where the absolute values of the roots of $b_f(s)$ are strictly smaller than $\al_1{+}1$, see \cite{sem}. A missing spectral number (counted with multiplicity) can be a root of $b_f(s)$ up to sign in the {\it semi-homogeneous\one} case, for instance if $f\eq x^8\pl y^8\pl x^5y^4$.
\end{rem}

\section{Proof of Theorem~\ref{T1}} \label{S5}
Let $e$ be the order of the semisimple part of the monodromy of $f$. The $V$-filtration on $A\eq\C\{x\}/(\dd f)$ can be indexed by $\Z$ identified with $\tfrac{1}{e}\,\Z\sst\Q$. Setting $F^k\!A\defs V^{k/e}\!A$, consider the mapping cone $(K^{\ssb},F^{\ssb})$ of the filtered morphism
\begin{equation*}
[f]:(A,F[-e])\to(A,F),
\end{equation*}
where $[f]F^p\!A\sst F^{p+e}\!A$. (This $F^{\ssb}$ is {\it not\one} the Hodge filtration.) We have decreasing filtrations $G,G'$ on $A$ such that
\begin{equation*}
G^p\!A\defs [f]F^p\!A,\q G'{}^p\!A\defs(F^p\!A\col[f])\eq{\rm Ker}([f]\col A\tos A/F^p\!A).
\end{equation*}
Note that $\bigcup_pG^p\!A\eq[f]A$. There is a spectral sequence $\{E_r^{p,q},\ddd_r\}_{r\ges0}$ associated with the filtered complex $(K^{\ssb},F^{\ssb})$ so that
\begin{equation} \label{5.1}
\aligned E_r^{p,-p}&=Z_r^{p,-p}/(B_r^{p,-p}\pl Z_{r-1}^{p+1,-p-1})\\&=F^p\!A/\bl(([f]F^{p-r-e+1}\!A\cap F^pA)+F^{p+1}\!A\br)\\&=\Gr_F^pA/G^{p-r-e+1}\Gr_F^pA,\endaligned
\end{equation}
\begin{equation} \label{5.2}
\aligned E_r^{p-r,-p+r-1}&=Z_r^{p-r,-p+r-1}/(B_r^{p-r,-p+r-1}\pl Z_{r-1}^{p-r+1,-p+r-2})\\&=\bl((F^p\!A\col[f])\cap F^{p-r-e}\!A\br)/\bl((F^p\!A\col[f])\cap F^{p-r-e+1}\!A\br)\\&=G'{}^p\Gr_F^{p-r-e}\!A,\endaligned
\end{equation}
see for instance \cite{Go}. Here we have by definition
\begin{equation*}
\aligned Z_r^{p,q}&\eq{\rm Ker}(\ddd\col F^pK^{p+q}\tos K^{p+q+1}/F^{p+r}K^{p+q+1}),\\
B_r^{p,q}&\eq{\rm Im}(\ddd\col F^{p-r+1}K^{p+q-1}\tos K^{p+q})\cap F^pK^{p+q},\endaligned
\end{equation*}
so that
\begin{equation*}
\aligned Z_r^{p,-p}&\eq F^p\!A,\q B_r^{p,-p}\eq[f]F^{p-r-e+1}\!A\cap F^pA,\\ Z_r^{p-r,-p+r-1}&\eq {\rm Ker}([f]\col F^{p-r-e}\!A\to A/F^p\!A),\q B_r^{p-r,-p+r-1}\eq 0.\endaligned
\end{equation*}
\sk
Since the differential $\ddd^j\col K^j\tos K^{j+1}$ vanishes for $j\nes{-}1$, there is an increasing filtration $U'$ on $\Gr_F^pK^0\eq\Gr_F^pA$ together with a decreasing filtration $U$ on $\Gr_F^pK^{-1}\eq\Gr_F^{p-e}A$ such that
\begin{align}
\Gr_F^pK^0/U'_{r-1}\Gr_F^pK^0\eq E_r^{p,-p},&\q U^r\Gr_F^pK^{-1}\eq E_r^{p,-p-1},
\label{5.3}\\
\Gr^{U'}_r\Gr_F^pK^0\eq{\rm Im}\,\ddd_r^{p-r,p+r-1},&\q\Gr_U^r\Gr_F^pK^{-1}\eq{\rm Coim}\,\ddd_r^{p,p-1},
\label{5.4}
\end{align}
where $U'_r\Gr_F^pK^0=0$ for $r\slt0$ and $U^r\Gr_F^pK^{-1}=\Gr_F^pK^{-1}$ for $r\less0$.
\sk
Comparing \eqref{5.3} with \eqref{5.1}--\eqref{5.2}, we see that
\begin{equation} \label{5.5}
U'_r\Gr_F^pK^0\eq G^{p-r-e}\Gr_F^pK^0,\q U^r\Gr_F^pK^{-1}=G'{}^{p+r}\Gr_F^pK^{-1}.
\end{equation}
Combined with \eqref{5.4}, these imply that
\begin{equation} \label{5.6}
\Gr_G^{p-r-e}\Gr_F^pA=\Gr_{G'}^p\Gr_F^{p-r-e}A=(\Gr_G^{ne-p}\Gr_F^{ne-p+r+e}A)^{\vee},
\end{equation}
since $G'{}^{\ssb}$ is the dual filtration of $G^{\ssb}[ne]$ with respect to $\Gr^F\Sb$ (see Remark\,\,\ref{R2.1}). This completes the proof of Theorem\,\,\ref{T1}.

\begin{rem}\label{R5.1}
Theorem\,\,\ref{T1} is very closely related to the {\it self-duality\one} of spectral sequence (see \cite[(1.3.1)]{De}, \cite[5.2.19–20]{mhp}), and may follow it if one can verify the compatibility of the duality isomorphisms with various isomorphisms of the spectral sequence, however we have {\it non-strict\one} two filtrations (the image filtration $V$ and its dual coimage filtration $I$), and the argument does not seem quite trivial. This self-duality does not seem to have been used for {\it never-degenerating\one} spectral sequences. They say that it seems to be used quite implicitly for the weight filtration in twistor theory, but this seems to be an $E_2$-degeneration case. Some people would expect that twistor theory can be applied to this setting as in the case of a computation of Kontsevich's complex; however, controlling the {\it non-strict\one} part of a {\it never-degenerating\one} spectral sequence seems very difficult. (Here one would use a not necessarily ``integrable" twistor modules, since they say that there is a still unclarified problem of ``half Tate twists" in the ``integrable" case, which are used in the computation of nearby cycles for instance.) If one considers the graded pieces of the pole order filtration $P^{\ssb}$ on the meromorphic de Rham complex, one can obtain a complex quasi-isomorphic to $K^{\ssb}$, but it does not seem easy to define the filtration $V$ on it (except the Newton nondegenerate case using the logarithmic complex, although the control of the non-strict part seems very hard). It seems quite unclear what kind of further structure can be expected.
\end{rem}

\begin{rem}\label{R5.2}
The estimate in Corollary\,\,\ref{C2} seems asymptotically optimal, considering the following families of polynomials
\begin{equation*}
f_{n,a}\defs\msum_{i=1}^n\,x_i^{an-1}\pl x_1^a\cdots x_n^a\,\,\,\,(a\,{\gg}\,0)\,\,\,\,\h{for}\,\,\,n\gess2.
\end{equation*}
Here
\begin{equation*}
\mu{-}\tau\eq\msum_{k=0}^{n-2}\,\tfrac{n!}{(n{-}k)!}\bl((n{-}k{-1})a\mi2\br)^{n-k},
\end{equation*}
and its $k\one $th summand corresponds to the graded piece
\begin{equation*}
\Gr^{\Wt}_{\ga+n}\Gr_V^{\ssb}{\rm Im}\,[f]\q\h{with}\q\ga\eq1{+}\tfrac{k+1}{e},\,\,e\eq an{-}1.
\end{equation*}
Indeed, setting
\begin{equation*}
\Si_k\defs\{{\bf p}\eq(p_1,\dots,p_k)\ins[1,n]^k\mid p_i\nes p_j\,\,\,\h{for any}\,\,i\nes j\},
\end{equation*}
the above graded piece is the direct sum of the vector spaces $V_{k,{\bf p}}$ (${\bf p}\ins\Si_k$) spanned by the monomials $x^{\nu}$ satisfying the following conditions:
\begin{equation*}
\aligned&\nu_j\eq ia{-}1\,\,\,\,\,\h{if}\,\,\,j\eq p_i\,\,(\exists\,i\ins[1,k]),\\ &\nu_j\ins[(k{+}1)a,na{-}3]\,\,\,\,\,\h{otherwise}.\endaligned
\end{equation*}
Note that the spectral number corresponding to $x^{\nu}$ is $(|\nu|{+}n)/e$.
\sk
Related to Corollary\,\,\ref{C2}, there are inequalities $\tfrac{\mu{-}\tau}{\mu}\slt\tfrac{1}{4}$ for $n\eq2$ (which was conjectured in \cite{DG}) and $\tfrac{\mu{-}\tau}{\mu}\slt\tfrac{1}{3}$ for $n\eq3$ assuming Durfee's conjecture \cite{Du} (which predicts the inequality $6p_g\less\mu$) for the latter, see \cite{Al0}. There are families of hypersurfaces showing that those two inequalities are asymptotically optimal, see \cite[Example 4.7]{Wa} for the case $n\eq3$, where the family is given by sufficiently general members of certain hypersurfaces (so the argument seems highly nontrivial). For $n\eq2$ however the above family is sufficient. Note that the limit of $\tfrac{\mu{-}\tau}{\mu}$ under $a\tos{+}\infty$ for the above family with $n\eq3$ is equal to $\tfrac{8}{27}$, which is strictly smaller than $\tfrac{1}{3}$.
\end{rem}

\section{Proof of Theorem~\ref{T2}} \label{S6}
The first assertion follows from Corollary\,\,\ref{C1.1}.
We see that $[f]\ins\Gr_I^{\al_1}\Gr_V^{\al^C_1}\!A$, and this does not vanish. Indeed, we have $I^{>\al_1}\Gr_V^{\al^C_1}\!A\eq0$ by the last assertion of Corollary\,\,\ref{C1.1}. This implies the inclusion
\begin{equation*}
\Gr_I^{\al_1}\Gr_V^{\al^C_1}\!A=I^{\al_1}\Gr_V^{\al^C_1}\!A\into\Gr_V^{\al^C_1}\!A,
\end{equation*}
where $[f]\ins\Gr_V^{\al^C_1}\!A$ does not vanish by the definition of $\al^C_1$.
So the last assertion of Theorem~\ref{T2} follows by applying Theorem\,\,\ref{T1} and using the equality $\al\mi\ga\eq n\mi\be$ for $\al\eq\al^C_1$, $\al\mi\ga\eq\al_1$, since $\al_1\eq n\mi\al_{\mu}$ by the symmetry of spectral numbers, see \eqref{3.4}.
(Note that Proposition\,\,\ref{P4.1} is not used in the proof of Theorem\,\,\ref{T2}, since Theorem\,\,\ref{T1} is sufficient. It does not seem very difficult to foresee a proof of the last assertion of Theorem~\ref{T2} if one looks at the proof of that proposition very carefully.)
This finishes the proof of Theorem\,\,\ref{T2}.

\section{Proof of Theorem~\ref{T4}} \label{S7}
It is sufficient to show the equality \eqref{3} as explained just after Theorem~\ref{T4}. We first prove the inequality $\si_{\!f}\less\si'_{\!f}$. By decreasing induction on $\be\sgt\al\mi\si_{\!f}$ (with $\al$ fixed) we can verify the vanishing
\begin{equation} \label{7.1}
\Gr_V^{\al}I^{\be}\!A\eq I^{\be}\Gr_V^{\al}A\eq\Gr_I^{\be}\Gr_V^{\al}A\eq0\q\h{if}\,\,\,\al\mi\be\slt\si_{\!f}.
\end{equation}
By increasing induction on $\al\less\be\pl\si_{\!f}$ (with $\be$ fixed) this implies that
\begin{equation} \label{7.2}
I^{\be}\!A\eq V^{\al}I^{\be}\!A,\,\,\h{that is,}\,\,\,I^{\be}\!A\sst V^{\al}\!A,\q\h{if}\,\,\,\al\mi\be\less\si_{\!f}.
\end{equation}
So the inequality $\si_{\!f}\less\si'_{\!f}$ follows. Assume the strict inequality holds here. By the definition of $\si_{\!f}$, there is $v\ins V^{\be}\!A$ with $[f]v\,{\notin}\,V^{>\al}\!A$ for some $\be\in\Q$ with $\al\defs\be\pl\si_{\!f}$ (since $I^{\be}\Gr_V^{\al}\!A$ is the image of $I^{\be}V^{\al}\!A$ in $\Gr_V^{\al}\!A$). However, this contradicts the strict inequality. So the equality \eqref{3} holds. This completes the proof of Theorem\,\,\ref{T4}.

\section{Proof of Proposition~\ref{P1}} \label{S8}
Set
\begin{equation*}
V_{\mhs N}^{\prime\one\al}\C\{x\}\defs V_{\mhs N}^{\al+\al_1}\C\{x\},\q V^{\prime\one\al}\!A\defs V^{\al+\al_1}\mhs A\q(\forall\,\al\ins\Q).
\end{equation*}
This is compatible with Remark~\ref{R3.4} in the case the inverse of each weight is an integer (for instance $f$ is a Brieskorn-Pham polynomial), since we have a {\it unique\one} maximal-dimensional compact face. We have the {\it strictly surjective canonical morphism of filtered rings}
\begin{equation} \label{8.1}
(\C\{x\},V_{\mhs N}^{\prime\ssb})\onto(A,V^{\prime\ssb}).
\end{equation}
Indeed, it is strictly surjective and the source is a filtered ring, hence so is the target. Recall that the ideal $V_{\mhs N}^{\al}\C\{x\}$ can be defined by the condition $\ell(\nu)\pl\al_1\gess\al$ for $x^{\nu}\ins\C\{x\}$ ($\nu\ins\N^n$) with $\ell(\nu)\defs\msum_{i=1}^n\,w_i\nu_i$ in the semi-weighted-homogeneous case, see Remark\,\,\ref{R3.5}. The assertions then follow from \eqref{8.1}. (This argument cannot be extended to the Newton non-degenerate case with several maximal-dimensional compact faces, since there is a big difference between $\pl{\bf 1}$ in \eqref{3.6} and ${+}\one\al_1$ in \eqref{3.9}.) The last assertion in the convenient Newton non-degenerate case follows from the remarks noted after \eqref{3.6}. This terminates the proof of Proposition\,\,\ref{P1}.

\section{Explicit calculations} \label{S9}
It is very important to verify abstract theories by explicit computations as much as possible, since any unexpected error might be hiding there.
We present a code in Singular verifying the graded symmetry of the missing spectral numbers of the $e$-Tjurina spectrum $\Spf^{(e)}(t)$ with $e\eq1$ or $2$ in the convenient Newton non-degenerate three-variable case as below. One notes the {\it vertices\one} of the Newton polyhedron together with a {\it higher\one} term $g$ {\it in the third line,} and $e$ (which must be either 1 or 2) is given {\it at the end of the second line.} The reader can {\it change these inputs\one} as long as they have {\it similar complexity.} The following code computes first the defining equations of the two-dimensional compact faces.
\ms
\vbox{\fontsize{8pt}{4mm}\sf\pv@LIB "linalg.lib"; LIB "sing.lib"; ring R=0, (x,y,z), ds; poly f,g,h,X,Z,Av,Sm,Vt;@
\pv@int a,b,c,e,i,j,k,m,n,v,p,q,r,s,t,Ms,Md,Mn,gn,dn,fl,lc,tj,mu,w; list M,T,spp; ideal I,J,K; e=1;@
\pv@intvec Inp=3,3,3, 12,0,0, 0,12,0, 0,0,12; a=7; b=2; c=2; g=x^a*y^b*z^c; /* a>=b>=c>1,a+b+2c=13 */@
\pv@v=size(Inp) div 3; intmat In[v][3]; for(i=1; i<=v; i++){for(j=1; j<=3; j++){In[i,j]=Inp[3*i+j-3];@
\pv@if(In[i,j]>m){m=In[i,j];}}} Ms=3*v; Mn=4*m; matrix vars[1][3]=x,y,z; matrix Co[3][1]; matrix@
\pv@Inv[3][3]; vector Id=[1,1,1]; matrix Y[v][1]; matrix AA[3][3]; matrix Eq[Ms][3]; fl=0; for(i=1;@
\pv@i<=v; i++){if(fl==0){sprintf("Input points:"); fl=1;} sprintf(" V
\pv@i,In[i,1],In[i,2],In[i,3]);} s=0; for(i=1; i<=v-2; i++){for(j=i+1; j<=v-1; j++){for(k=j+1; k<=v;@
\pv@k++){for(p=1;p<=3;p++){AA[1,p]=In[i,p]; AA[2,p]=In[j,p]; AA[3,p]=In[k,p];} if(det(AA)!=0){Co=@
\pv@inverse(AA)*Id; Y=In*Co; for(q=1; q<=v && sum(Y,q)>=1;q++){;} if(q>v){X=sum(vars*Co); for(r=1;@
\pv@r<=s & X!=Eq[r,1]; r++){;} if(r>s){s++; Eq[s,1]=X; Eq[s,2]=cleardenom(X); Eq[s,3]=int(number(@
\pv@diff(Eq[s,2],x))/number(diff(Eq[s,1],x))); if(Eq[s,3]>Md){Md=int(Eq[s,3]);} if(s==1){printf(@
\pv@"2-dim compact faces:");}sprintf("S
\pv@Eq[s,3]); for(p=1; p<=v; p++){if(subst(Eq[s,1],x,In[p,1],y,In[p,2],z,In[p,3])==1){sprintf(@
\pv@" V
\msn
We can get the Steenbrink and the $e$-Tjurina spectral numbers by computing the Newton filtration, and verify the symmetry of the Steenbrink spectrum. (If we apply this to a {\it Newton degenerate\one} case, for instance if $f\eq xyz(x{+}y{+}z)\pl x^5\pl y^5\pl z^5$, it is mathematically incorrect and the computation by Singular may become unreliable.)
\ms
\vbox{\fontsize{8pt}{4mm}\sf\pv@lc=1; for(i=1; i<=s; i++){lc=lcm(lc,int(Eq[i,3]));} matrix G[2][s*Md]; for(i=1;i<=s;i++)@
\pv@{for(j=1; j<Eq[i,3]; j++){G[1,(i-1)*Md+j]=u*j/Eq[i,3];}} fl=0; for(i=1; i<=s*Md&&fl==0;i++)@
\pv@{X=1; if(i>1){Z=G[2,i-1];} else {Z=0;} for(j=1; j<=s*Md; j++){if(G[1,j]<X && G[1,j]>Z){X=G[1,j];@
\pv@}} if(X==1){fl=1;} else {G[2,i]=X;}} gn=i-1; G[2,gn]=1; matrix Gt[1][3*gn]; matrix Di[3*gn][2];@
\pv@matrix S[3*gn][4]; intmat FD[3*gn][3*gn]; intmat GD[3*gn][3*gn]; intmat C[gn*gn][6]; for(i=0; i<3;@
\pv@i++){for(j=1; j<=gn; j++){Gt[1,gn*i+j]=i+G[2,j];}} h=0; for(i=1; i<=v; i++){h=h+i*x^(In[i,1])*y^(@
\pv@In[i,2])*z^(In[i,3]);} f=h+g; tj=tjurina(f); mu=milnor(f); matrix TSp[1][mu]; matrix Sps[mu][2];@
\pv@if(mu!=milnor(h)){printf(" g error!!!"); exit;} M[3*gn]=jacob(f); I=(jacob(f),f^e); T[3*gn]=I;@
\pv@for(p=3*gn-1; p>0; p--){M[p]=M[p+1]; T[p]=T[p+1]; for(q=1; q<=s; q++){if(denominator(number(Eq[@
\pv@q,3]*Gt[1,p]))==1){for(i=1; i<=Mn; i++){for(j=1; j<=Mn; j++){m=int(Gt[1,p]*Eq[q,3]-i*diff(Eq[q,2],@
\pv@x)-j*diff(Eq[q,2],y)); n=int(diff(Eq[q,2],z)); if(m
\pv@fl==0; r++){if(subst(Eq[r,1]-Eq[q,1],x,i,y,j,z,k)<0){fl=1;}} if(fl==0){M[p]=(M[p],x^(i-1)*y^(j-1)*@
\pv@z^(k-1)); T[p]=(T[p],x^(i-1)*y^(j-1)*z^(k-1));}}}}} M[p]=std(M[p]); T[p]=std(T[p]);}} for(p=3*gn;@
\pv@p>0; p--){Di[p,1]=vdim(std(M[p])); Di[p,2]=vdim(std(T[p])); if(p<3*gn){S[p,1]=Di[p+1,1]-Di[p,1];@
\pv@S[p,2]=Di[p+1,2]-Di[p,2]; S[p,3]=S[p,1]-S[p,2];}} t=0; q=0; for(p=1;p<=3*gn-1; p++){if(S[p,1]>0)@
\pv@{if(t==0){printf("Spec. numb. with mult. (also for T-spec. and C-spec.)"); t=@
\pv@1;} sprintf(" 
\pv@Sps[q,2]=S[p,1];}} for(i=1; 2*i<=q; i++){if(Sps[i,1]+Sps[q+1-i,1]!=3 || Sps[i,2]!=Sps[q+1-i,2])@
\pv@{sprintf(" Symmetry of spectrum fails!"); exit;}} sprintf(" Symmetry of spectrum holds."); w=q;@}
\msn
We verify the graded symmetry of missing spectral numbers for the $e$-Tjurina spectrum by calculating the spectral sequence. This might take extremely long even if $f$ does not seem so complicated; in this case it is better to stop the computation and change the input. Note however that one can {\it never\one} determine the shifts of degrees by looking only at the missing spectral numbers even if one can foresee them.
\ms
\vbox{\fontsize{8pt}{4mm}\sf\pv@r=0; for(p=1;p<3*gn; p++){fl=0; for(q=p-1; q>0 && fl==0; q--){if(Gt[1,q]+1<=Gt[1,p] && S[p,3]>0@
\pv@&& S[q,1]>0){I=(f^e*M[q],M[p+1]); J=intersect(I,M[p]); FD[p,q]=vdim(std(M[p+1]))-vdim(std(J));@
\pv@for(k=q+1; k<3*gn && FD[p,k]==0; k++){;} GD[p,q]=FD[p,q]-FD[p,k];@
\pv@if(Gt[1,p]-Gt[1,q]>=1 && GD[p,q]!=0){r++; C[r,1]=int(Gt[1,p]*lc); C[r,2]=GD[p,q];@
\pv@C[r,3]=int(lc*(Gt[1,p]-Gt[1,q]));} if(FD[p,q]==S[p,3]){fl=1;}}}} fl=0; for(i=1; i<=r && fl==0;@
\pv@i++){m=3*lc; if(i>1){n=C[i-1,4];} else {n=0;} for(j=1; j<=r; j++) {if(C[j,3]<m && C[j,3]>n){m=@
\pv@C[j,3];}} if(m==3*lc){fl=1;} else {C[i,4]=m;}} dn=i-2; m=0; t=0; for(i=1; i<=dn; i++){fl=0; p=0;@
\pv@for(j=1; j<=r; j++){if(C[j,3]==C[i,4]){if(fl==0){if(t==0){printf("Shift deg of miss. sp. numb."@
\pv@); t=1;}sprintf("Shifted by 
\pv@lc,C[j,2]); p++; C[p,5]=C[j,1]; C[p,6]=C[j,2];}} for(q=1; q<=p div 2; q++){if(C[q,5]+C[p+1-q,5]@
\pv@!=C[i,4]+3*lc || C[q,6]!=C[p+1-q,6]){sprintf(" Graded symmetry error!"); m++;}}} if(m==0){@
\pv@sprintf(" Graded symmetry holds.");};@}
\msn
Here we see many examples of missing spectral numbers $\al_i,\al_j$ such that $\al_i\eq\al_j$ but their corresponding $\ga$ differ.
Concerning the non-simplicial Newton polyhedron case, one may get a Newton {\it degenerate\one} polynomial, although this does not occur usually in the case where the symmetry of the spectrum holds. This problem may be solved by {\it changing the order of vertices.} One can see $h$ by typing ``{\smaller\sf\verb@h;@}".

\end{document}